\documentclass[11pt,a4paper]{article}
\usepackage{graphics}
\usepackage{graphicx}
\usepackage{amsmath}
\usepackage{amssymb}
\usepackage{authblk}
\usepackage[numbers,sort&compress]{natbib}

\begin{document}
\title{Effects of time-delay in a model of intra- and inter-personal motor coordination}

\author[1]{Piotr S\l{}owi\'{n}ski}
\author[,1]{Krasimira~Tsaneva-Atanasova\thanks{K.Tsaneva-Atanasova@exeter.ac.uk}} 
\author[2]{Bernd~Krauskopf}
\affil[1]{Department of Mathematics,
        College of Engineering, Mathematics and Physical Sciences,
	University of Exeter,
	Exeter, Devon, EX4 4QF, UK} 
\affil[2]{Department of Mathematics,
	University of Auckland, Auckland 1142, New Zealand}
\maketitle

\abstract{Motor coordination is an important feature of intra- and inter-personal interactions, and several scenarios --- from finger tapping to human-computer interfaces --- have been investigated experimentally. In the 1980, Haken, Kelso and Bunz formulated a coupled nonlinear two-oscillator model, which has been shown to describe many observed ascpects of coordination tasks. We present here a bifurcation study of this model, where we consider a delay in the coupling. The delay is shown to have a significant effect on the observed dynamics. In particular, we find a much larger degree of bistablility between in-phase and anti-phase oscillations in the presence of a frequency detuning.} 

\section{Introduction}
\label{intro}

Many joint-action tasks demand some degree of movement coordination. Moreover, the degree of movement coordination plays an important role in inter-personal interactions; e.g., it can affect the level of affiliation between interacting people \cite{Oullier2008}. In the case of (near) periodic movements the collective patterns of coordination are well captured by the properties of the relative phase, $\phi_\mathrm{rel}$, between the individual coupled oscillating subsystems \cite{Kelso1997}. In the case of two coupled oscillators the simplest coordination pattern is observed when the phase of the two oscillators coincide to give in-phase monostable coordination (where $\phi_\mathrm{rel}=0$). Monostable anti-phase coordination (where $\phi_\mathrm{rel}=\pi$) can also occur, and an example of such monostable behaviour is observed in competitive games \cite{Bourbousson2010,Duarte2012}. In many real systems stable anti-phase coordination coexists with stable in-phase coordination \cite{Haken1985,Kelso1997,Mortl2012,Warren2006}. 
The development of the well-known Haken-Kelso-Bunz (HKB) \cite{Haken1985} model (see \eqref{eqn4} below) has been inspired by the in-phase and anti-phase coordination dynamics observed in bimanual coordination experiment \cite{Kelso1983}.
In the past 30 years the HBK model has been widely applied as a paradigm for studying dynamics of intra- and inter-personal motor coordination. It was found to be representative of a wide range of human movement experiments \cite{Calvin2011}, suggesting that the dynamics observed in the HBK model are somehow fundamental  \cite{Kelso1997}.

The importance of perceptual-motor delays have been acknowledged in the original HKB model paper \cite{Haken1985} and later discussed in \cite{Peper1998}. However, the possible role of time delays in the HKB model has not been analysed in a systematic way. More generally, there have been very few studies concerning models of motor coordination that incorporate time delays. A study of an excitator model with time delay has been presented in \cite{Banerjee2007}, and a model of relative phase dynamics with time delay can be found in \cite{Tass1995,Tass1996}. More recently, simulations of the HKB model with delays have been used to interpret interpersonal coordination patterns observed in the interaction differences between schizophrenic patients and healthy participants \cite{Varlet2012}. 

The full HKB model \cite{Haken1985} with delays in the coupling function can be written as:
\begin{align}
  \label{eqn4}
  \dot{x_1}(t) = & y_1(t), \\
  \nonumber
  \dot{x_2}(t) = & y_2(t), \\
  \nonumber
  \dot{y_1}(t) = &- \left(y_1(t) \left( \alpha x_1(t)^2  + \beta  y_1(t)^2 - \gamma \right)+\omega ^2 x_1(t)  \right)\\
  \nonumber
  & + \left(a_1+b_1\left(x_1(t)-x_2(t-\tau_1)\right)^2\right)\left(y_1(t)-y_2(t-\tau_1)\right), \\  \nonumber
  \dot{y_2}(t) = &- \left(y_2(t) \left( \alpha x_2(t)^2  + \beta  y_2(t)^2 - \gamma \right)+(\omega+\Delta)^2 x_2(t)  \right) \\  \nonumber
 & + \left(a_2+b_2\left(x_2(t)-x_1(t-\tau_2)\right)^2\right)\left(y_2(t)-y_1(t-\tau_2)\right).
\end{align}
In this system of four first-order differential equations with delays the variables $x_1(t)$ and $x_2(t)$ represent the positions and $y_1(t)$ and $y_2(t)$ the velocities of the two individual oscillators at time $t$, where time is measured in seconds. 
The parameters $\tau_1$ and $\tau_2$ denote time delays arising due to cognitive processes and/or physiological properties of the neuromuscular system. Hence, $x_2(t-\tau_1)$, $y_2(t-\tau_1)$ and $x_1(t-\tau_2)$, $y_1(t-\tau_2)$ represent the positions and velocities at times $t-\tau_1$ and $t-\tau_2$, respectively. Further, $\omega \in \mathbb{R}^+$ is the frequency of the oscillations (either the natural or eigenfrequency or the external pacing) and $\Delta$ is the
detuning of the second oscillator with respect to $\omega$. The parameters $\alpha, \beta, \gamma \in \mathbb{R}$ in \eqref{eqn4} govern the intrinsic dynamics of the single HKB oscillator; these include the Rayleigh oscillator term $\beta y(t)^3$, as well as linear and  Van der Pol-type nonlinear damping terms $\gamma y(t)$ and $\alpha x(t)^2 y(t)$, respectivley. Finally, the parameters $a_1, a_2, b_1$ and $b_2$ control the strength of coupling between the two HBK oscillators. Indeed, for $\tau_1 = \tau_2 = 0$ system \eqref{eqn4} is exactly the four-dimensional HBK ODE. On the other hand, the HBK model with $\tau_1, \tau_2 \neq 0$ is a delay differential equation (DDE); hence, it has as its phase space the infinite-dimensional space of continuous functions from the (maximal) delay interval into the four-dimensional $(x_1,x_2,y_1,y_2)$-space \cite{Diekmann1995}.

In this paper, we present a bifurcation study, by means of numerical continuation, of the HKB model with time delays that investigates further and explains simulation results in \cite{Varlet2012}. Since delay differential equations can exhibit richer dynamics than ordinary differential equations, this type of study constitutes an important step towards a deeper understanding of more realistic models for motor coordination. More specifically, we find that time delay has a large effect with regard to regions of multistability between in-phase and anti-phase solutions. Our results provide a better understanding of the experiments and, importantly, allow us to make experimentally testable predictions. All computations were performed with the latest version of the continuation package DDE-Biftool v3.1 \cite{DDEBIFTOOL} under Matlab.

More specifically, we consider \eqref{eqn4} with equal coupling strengths $a=a_1=a_2$ and  $b=b_1=b_2$. This is usually regarded as a model of intra-personal coordination, such as coordination of the left and right hand of the same person; therefore, we also assume an equal delay $\tau=\tau_1=\tau_2$. Throughout, the intrinsic parameters of the two oscillators are fixed to values estimated directly from experimental data for wrist movements \cite{Kay1987}, namely to $\alpha=12.457$, $\beta=0.007905$ and $\gamma=0.641$; these values were also used in a relatively recent study of a virtual partner interaction \cite{Kelso2009}.

\section{Influence of the pacing frequency}
\label{sec:1}

We first investigate how the solutions of \eqref{eqn4} change with increasing frequency $\omega$ of the oscillations, which is commonly interpreted as an increase in pacing frequency during an experiment. Here we assume that there is no detuning and, hence, set $\Delta = 0$. To evaluate the influence of the delay $\tau$, we compare one-parameter bifurcation diagrams in $\omega$ of \eqref{eqn4} for $\tau=0$ and for $\tau=0.14$. The latter value of the time delay was chosen because it is well within the range $70-150\,$ms that has been reported for responses to continuous stimuli \cite{Loram2006}; it is also close to the value of $170\,$ms that was measured for refractory reactions \cite{Miall1993}. 

\begin{figure}
\centering
\includegraphics{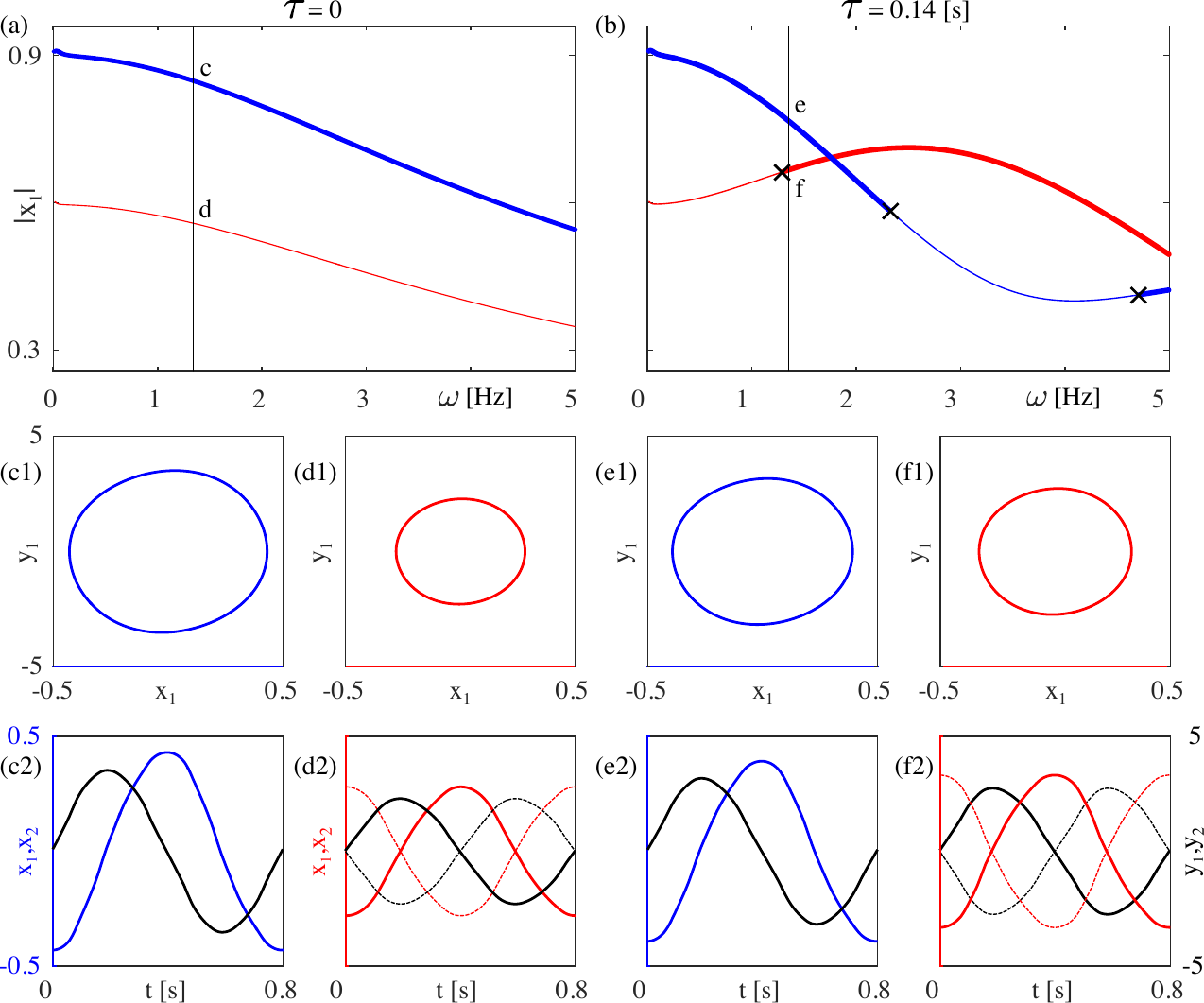}
\caption{One-parameter bifurcation diagrams in $\omega$ of \eqref{eqn4} for $\tau=0$ (a) and for $\tau=0.14$ (b), shown in terms of the amplitude $|x_1|$. Blue curves indicate in-phase solutions and red curves anti-phase solutions; thick curves denote stable and thin curves  unstable solutions, and black crosses (x) are branch points. Panels (c)--(f) show the periodic solutions for $\omega=1.3$ (grey vertical line), where (c1)--(f1) illustrate phase portraits and (c2)--(f2) are time series over a single period. Coloured curves in (c2)--(f2) indicate position $x$ (left axis) and black curves velocity $y$ (right axis); the thin dashed curves of the anti-phase solutions in (d2) and (f2) represent $x_2$ and $y_2$, respectively. Here $\alpha=12.457$, $\beta=0.007095$, $\gamma=0.641$, $a=-0.2$, $b=0.2$ and $\Delta = 0$.}
\label{fig:w_1D} 
\end{figure}

The two bifurcation diagrams are shown in Fig.~\ref{fig:w_1D} in terms of the amplitude $|x_1|$ of the position of the first oscillator; note that $|x_1|=|x_2|$ since the two coupled oscillators are identical. In-phase solutions are depicted as blue curves and anti-phase solutions as red curves, while their stability is indicated by the thickness of the respective curve. In Fig.~\ref{fig:w_1D} we set the coupling coefficients to $a=-0.2$ and $b=0.2$  as in \cite{Haken1985}. 

Figure~\ref{fig:w_1D}(a) shows that for $\tau=0$ the in-phase solution is the only stable solution over the entire $\omega$-range shown. Hence, without delay we do not find any bistability for the chosen experimentally validated values of $\alpha$, $\beta$ and $\gamma$. (We remark that these are very different from the choice $\alpha=0$, $\beta=1$ and $\gamma=1$  that was used when bistability was found in the original paper \cite[Fig.~7]{Haken1985}.) The corresponding bifurcation diagram for $\tau=0.14$ in Fig.~\ref{fig:w_1D}(b) demonstrates that the introduction of time delays may significantly affect the stability properties of solutions of the HKB model \eqref{eqn4}. There are now two regions of bistability, one for low values of $\omega\in(1.28, 2.33)$ and the other for high values of $\omega>4.7$. Stability of the periodic solutions is gained or lost via branch points (where a single real Floquet multiplier crosses through 1; note that there is always an additional trivial Floquet multiplier 1). We remark that the stability properties of the periodic solutions of \eqref{eqn4} depend also on the strength of coupling, but this is not explored further here. 

Panels (c)--(f) of Fig.~\ref{fig:w_1D} show phase portraits and time series of representative in-phase and anti-phase solutions for $\omega=1.3$ (indicated by the vertical line in the bifurcation diagrams in panels (a) and (b)). This value has been chosen because it lies within the range of frequencies at which bistability between in-phase and anti-phase solutions has been observed in experiments \cite{Kay1987,Molenaar2003}. Panel (c1) shows the stable in-phase periodic orbit for $\tau = 0$ in the $(x_1,y_1)$-plane, and panel (c2) shows its positions and velocities over one period. The unstable solution that exists for the same value of $\omega=1.3$ is shown in the same way in panels (d1) and (d2), and note that now the two positions and velocities are indeed exactly in anti-phase. Panels (e) and (f) also show the in-phase and anti-phase solutions at $\omega=1.3$ but now for $\tau=0.14$. These solutions are quite similar to those for $\tau = 0$, but they are both stable.

\section{Regions of multistability in the $(a,\tau)$-plane}
\label{sec:at}

\begin{figure}
\centering
\includegraphics{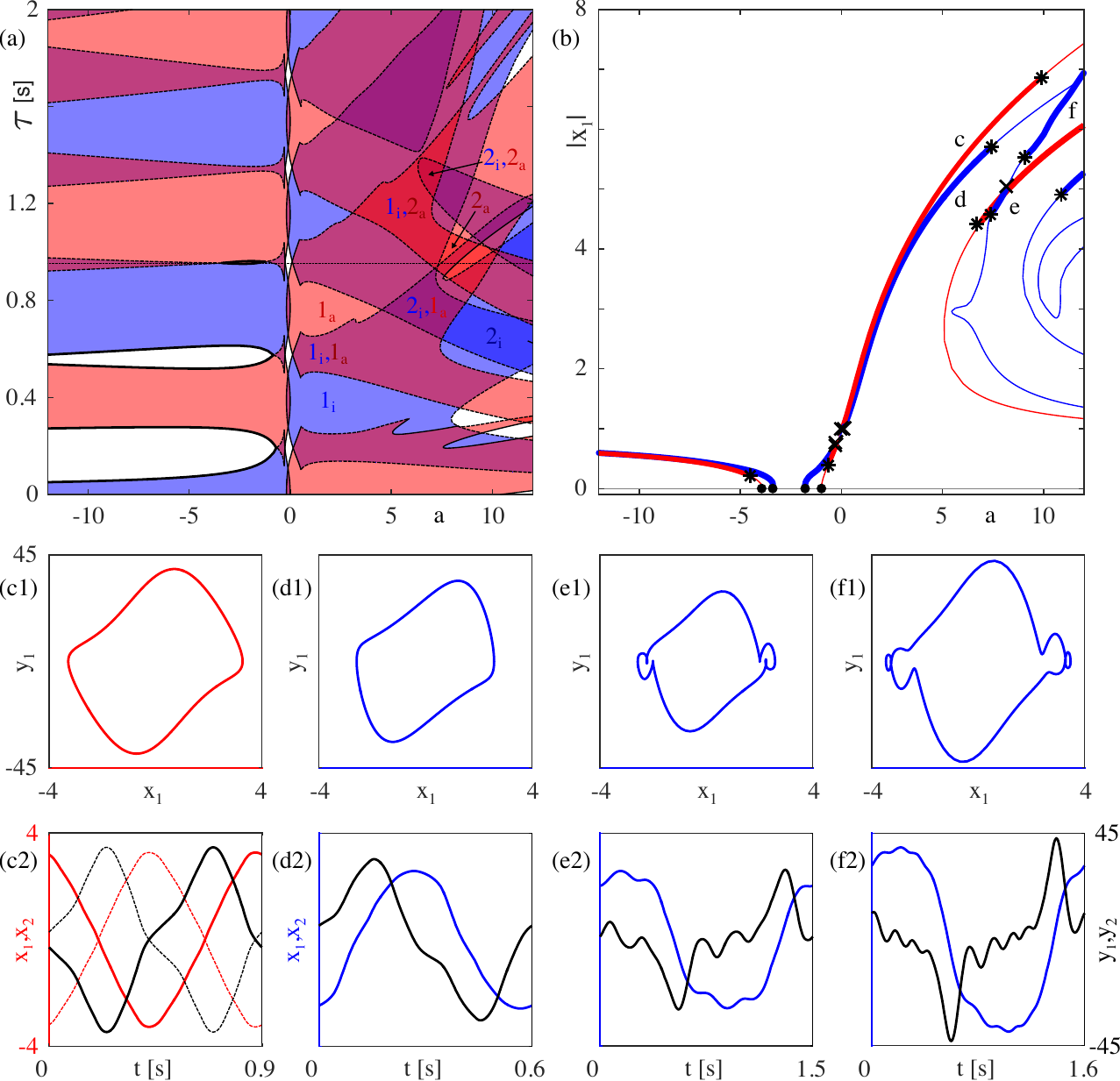}
\caption{Two-parameter bifurcation diagram of \eqref{eqn4} in the $(a,\tau)$-plane for $\omega=1.3$ (a), associated one-parameter bifurcation diagram in $a$ for $\tau=0.95$, and periodic solutions for $\tau=0.95$ and $a=5$ (c)--(d), $a=8$ (e) and $a=11.5$ (f). 
Blue regions in panel (a) indicate stable in-phase solutions and red regions stable anti-phase solutions; multistability is indicated by different shades of red, blue and purple. Thick black curves are loci of Hopf bifurcations, thin black curves are loci of branch points and dashed black curves are loci of torus bifurcations; in panel (b) these bifurcations are indicated by 
dots (\textbullet), crosses (x) and asterisks (*), respectively. Here $\alpha=12.457$, $\beta=0.007095$, $\gamma=0.641$, $b=1$ and $\Delta = 0$.}
\label{fig:at} 
\end{figure}

In order to gain insight into the effect of the time delay $\tau$ on the stability of the in-phase and anti-phase solutions we compute the two-parameter bifurcation diagram of \eqref{eqn4} in the $(a,\tau)$-plane of coupling coefficient and time delay. We investigate the dependence of the solutions on the coupling parameter $a$, where $b=1$ is fixed. This choice is based on the observation that the bistability between in-phase and anti-phase solutions in the HKB model without delay depends on $a$ for a range of values for $b$; see \cite{Avitabile2016}.
In this bifurcation analysis we fix the frequency to $\omega=1.3$, as in Fig.~\ref{fig:w_1D}(c)--(f), and explore a quite large range of $\tau$. Namely, we not only consider values of $\tau$ that are within the physiologically relevant range of time delays observed in human physiology and cognition $\tau<0.8$ \cite{Thorpe1996,FabreThorpe2001,Michel2004,Li2012}, but also larger time delays up to  $\tau=2$. The latter range is relevant in coordination games with virtual partners \cite{Kelso2009} or coordination between two people who are communicating (with lags) over long distances. Similarly, we explore a relatively large range of the coupling coefficient $a$ because there are no reliable estimates from experimental data \cite{Schoner1986,Post2000,Molenaar2003}. 

Figure~\ref{fig:at}(a) shows the two-parameter bifurcation diagram of \eqref{eqn4} in the $(a,\tau)$-plane for $\omega=1.3$, with regions of stable in-phase and anti-phase solutions coloured accordingly; labels indicate the number of stable in-phase (subscript i) and stable anti-phase (subscript a) solutions. In the white regions there are neither stable in-phase nor stable anti-phase solutions. The different regions are bounded by curves of Hopf bifurcations, branch points and torus bifurcations. It is important to keep in mind that we show only those  parts of the respective bifurcation curves that bound regions of stable solutions. Moreover, we emphasise that we are only considering here solutions that appear via Hopf bifurcations of the trivial steady state (0,0,0,0). Additional periodic solutions (for example, those associated with symmetry-breaking  or locking on tori) may exist in the the blue, red and white regions of the $(a,\tau)$-plane; however, these are beyond the scope of this study. 

In Fig.~\ref{fig:at}(a) there are Hopf bifurcations only for $a<0$, while we find that periodic solutions lose or gain stability at branch points and torus bifurcations. The curves of Hopf bifurcations associated with the in-phase and anti-phase solutions intersect at double Hopf points, which give rise to two branches of torus bifurcations. There are many parameter values where regions of stable in-phase and anti-phase solutions overlap, giving rise to bi- or multistability. Figure~\ref{fig:at}(a) clearly shows that, as $\tau$ increases, regions of stable in-phase and anti-phase solutions alternate. This can be understood intuitively by considering the relationship between the time delay $\tau$ and the  frequency $\omega$: stable in-phase or anti-phase can be found separated by a time shift equal to half of the period of the oscillations; for $\omega=1.3$ this corresponds to a distance in $\tau$ of $0.3846$.

To illustrate the multistability further, Fig.~\ref{fig:at}(b) presents the one-parameter bifurcation diagram of \eqref{eqn4} for increasing values of $a$ and fixed $\tau=0.954$, which is a value (grey horizontal line in panel (a)) for which there is a complicated structure of different periodic solutions. There are four branches of periodic solutions that appear via Hopf bifurcations of the steady state: two of them extend to the left and two to the right. We further found four isolated branches: three of in-phase and one of anti-phase periodic solutions. The three branches of in-phase solutions are actually connected by fold bifurcations at different values of $\tau$ with the branch of the in-phase solutions that emanates from a Hopf bifurcation and extends to the right; the same applies to the separate branch of the anti-phase solutions. 
Taken together, panels (a) and (b) of Fig.~\ref{fig:at} demonstrate that there is a large degree of multistability, which increases with $\tau$ and is much more pronounced for a positive coupling strength $a$. 

Figure~\ref{fig:at}(c)--(f) are examples of the different types of stable periodic solutions that can be found along the solution branches in Fig.~\ref{fig:at}(b); they are shown again as periodic orbits in the $(x_1,y_1)$-plane and time series of positions and velocities over one period. 
The periodic solutions in panels (c) and (d), along the solution branches emanating from Hopf bifurcations, are very similar to those from Fig.~\ref{fig:w_1D}(e) and (f). The periodic solutions in Fig.~\ref{fig:at}(e) and (f), on the other hand, show an
interesting new feature, namely an increasing number of smaller maxima, which correspond to little loops in projection onto the $(x_1,y_1)$-plane. We remark this phenomenon is not related to changes in the stability of the solutions; it is similar to spiking observed in model of coupled neural populations presented in \cite{Marten2009}. A detailed analysis of the emergence of maxima in \eqref{eqn4} is beyond the scope of this paper.

\section{Frequency detuning and relative phase}
\label{sec:awd}

\begin{figure}
\centering
\includegraphics{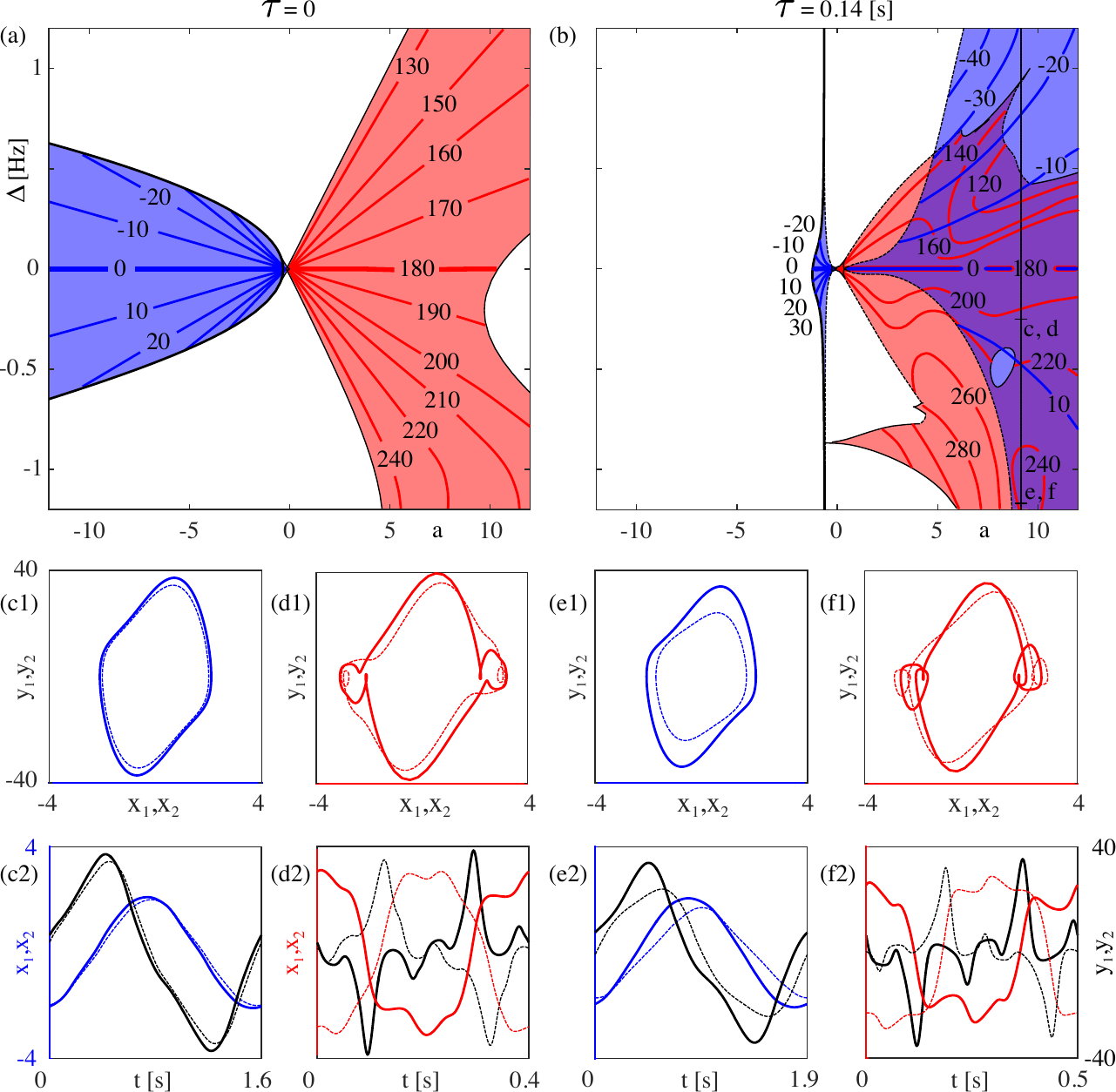}
\caption{Two-parameter bifurcation diagrams of \eqref{eqn4} in the $(a,\Delta)$-plane for $\tau=0$ (a) and for $\tau=0.14$ (b). Curves and regions are as in the previous figures; also shown are contours of relative phase $\phi_\mathrm{rel}$ (in degrees) of stable periodic solutions. Panels (c)--(f) show periodic solutions from the region of bistability in panel (b) for $a=9.1685$, namely for $\Delta=-0.25$ with $\phi_\mathrm{rel}=5.4^{\circ}$ (c) and with $\phi_\mathrm{rel}=205.2^{\circ}$ (d), and for $\Delta=-1.14$ with $\phi_\mathrm{rel}=17.9^{\circ}$ (e) and with $\phi_\mathrm{rel}=240.3^{\circ}$ (f). Here $\alpha=12.457$, $\beta=0.007095$, $\gamma=0.641$, $b=1$ and $\omega=1.3$.}
\label{fig:awd} 
\end{figure}

We now consider the case that there is a difference between the frequencies of the two oscillators, which is expressed in \eqref{eqn4} as the detuning $\Delta$ of the second oscillator with respect to the frequency $\omega$ of the first oscillator. A detuning between the oscillators is one of the main parameters that can be controlled in experiments. For example, it is a common experimental practice in studies of inter-personal coordination to detune the intrinsic frequencies between the pendula driven by the wrist movement by $|\Delta|<0.4$ Hz \cite{Schmidt1997,Varlet2012,Varlet2014}. Furthermore, even for the case of intra-personal coordination of wrists movements (or if the pendula in the experimental set-up are identical) one should expect a small detuning due to small differences between the left and the right hand (or due to individual differences between participants) \cite{Kay1987,Schmidt1997}. As a result of detuning ($\Delta \neq 0$) the relative phase $\phi_\mathrm{rel}$ between the two oscillators may take values other than 0 and $\pi$.
In other words, When $\Delta$ is increased from 0, the in-phase and anti-phase solutions generalise to solutions with a phase difference $\phi_\mathrm{rel}$ near 0 and near $\pi$, respectively; these intermediate-phase solutions are still referred to as (generalised) in-phase and anti-phase solutions for simplicity, provided the change in $\phi_\mathrm{rel}$ is reasonably small. To determine the phase of an oscillator we consider the Hilbert transform of a periodic orbit as a proto-phase, which we then transform into an observable-independent phase $\phi$
that grows linearly in time; see \cite{Kralemann2008} for more details of this technique. 

Figure~\ref{fig:awd} illustrates the two-parameter bifurcation diagram of \eqref{eqn4} in the $(a,\Delta)$-plane for $\tau =0$ and for $\tau = 0.14$; also shown are examples of periodic solutions for $\tau = 0.14$. In the bifurcation diagrams in panel (a) and (b) the central case $\Delta=0$ signifies that the two oscillators have the same frequency, as was the case in previous sections. The blue regions indicate stable (generalised) in-phase periodic solutions that originate from the in-phase oscillation for $\Delta=0$ with $\phi_\mathrm{rel}=0$, while the red regions indicate stable (generalised) anti-phase periodic solutions that originate from the anti-phase oscillation for $\Delta=0$ with $\phi_\mathrm{rel}=\pi$. The curves are contours of equal relative phase $\phi_\mathrm{rel}$, which is given and labelled in term of the phase angle. As before, in the white regions of Fig.~\ref{fig:awd}(a) and (b) there are no stable in-phase nor anti-phase solutions that appear via Hopf bifurcations of the trivial steady state (0,0,0,0). 

Figure~\ref{fig:awd}(a) for $\tau=0$ shows a region of stable periodic solutions with phases near $\phi_\mathrm{rel}=0^\circ$, which exists for $a<0$ and features $\phi_\mathrm{rel}$ up to about $\pm 40^\circ$, and a region of stable periodic solutions with phases near $\phi_\mathrm{rel}=180^\circ$, which exists for $a>0$ and features $\phi_\mathrm{rel}$ in the range $180^\circ \pm 50^\circ$. Notice that the contours of $\phi_\mathrm{rel}$ are very regular and practically symmetric with respect to $\Delta = 0$. The blue region is bounded mostly by loci of Hopf bifurcations and the red region is bounded entirely by loci of branch points.

Figure~\ref{fig:awd}(b) for $\tau=0.14$ shows that the introduction of delay has a significant effect on the two typs of stability regions in the $(a,\Delta)$-plane. First of all, the blue region of stable in-phase periodic solutions for $a<0$ almost disappears. On the other hand, there is now a large region of stable in-phase periodic solutions for $a>0$, which overlaps to a large extend with the red region of stable anti-phase periodic solutions. This results in a considerable and experimentally significant region of bistability between in-phase and anti-phase coordination regimes. Again, these regions are bounded by loci of branch points as well as Hopf and torus bifurcations. Notice further that the stability regions for $a>0$ and the associated contours of relative phase are no longer symmetric with respect to $\Delta = 0$. In particular, the variation in  $\phi_\mathrm{rel}$ of the stable anti-phase solutions near $\Delta=0$ is larger than that of the in-phase solutions. These results are consistent with experimental findings demonstrating that intra- as well as inter-personal anti-phase coordination is less stable than the in-phase coordination \cite{Kay1987,Schmidt1997,Varlet2012,Varlet2014}.

Finally, panels (c)--(f) of Fig.~\ref{fig:awd} are examples of periodic solutions for $\tau=0.14$ that can be found in panel (b) for $a=9.1685$ (indicated by the grey vertical line through the region of bistability). For a quite small frequency detuning of $\Delta= -0.25$ one finds the stable in-phase periodic solution with $\phi_\mathrm{rel}=5.4^{\circ}$ and the stable anti-phase periodic solution with $\phi_\mathrm{rel}=205.2^{\circ}$ shown in Fig.~\ref{fig:awd}(c) and (d), respectively.
For a larger detuning of $\Delta=-1.14$ one finds the quite similar stable periodic solutions with $\phi_\mathrm{rel}=17.9^{\circ}$ in (e) and with $\phi_\mathrm{rel}=240.3^{\circ}$ in (f).  This demonstrates that for sufficiently strong coupling it is possible to achieve coordination even in the face of a large frequency detuning between the two oscillators. Notice from the projections in (c1)--(f1) that the traces of the two oscillators in the $(x_1,y_1)$-plane and the $(x_2,y_2)$-plane are no longer identical, and compare with the time series in  (c1)--(f1).  Moreover, the stable periodic solutions originating from the anti-phase solutions for $\Delta = 0$ exhibit behaviour similar to the in-phase solutions for large $a$ in Fig. \ref{fig:at}(e)--(f). Interestingly, the additional maxima are much clearer in the time series of the velocities (black curves) than in the time series of positions. This means that they might be difficult to detect in experimental (and, hence, noisy) position data.

\section{Conclusions}

We presented a bifurcation study of the HKB model for physically relevant parameter settings that are consistent with experimental data, and with delay in the coupling. The focus was on the effect of the delay on stable in-phase and anti-phase periodic solutions, where we considered also the case when a frequency detuning is present. More specifically, the comparison of two-parameters bifurcation diagrams, without and with delay, clearly demonstrates that time delays may have a significant effect, especially on the size and features of regions of bistability between in-phase and anti-phase oscillations. Furthermore, our analysis of the HKB model with time delay and frequency detuning provides a possible dynamical explanation why the anti-phase solutions are less stable in practice. Namely, already small differences between the intrinsic frequencies of the two oscillators result in a much larger deviations from the desired coordination pattern for the anti-phase compared to the in-phase solutions.

Our study is merely a first attempt at a bifurcation analysis of the full four-dimensional HKB model with time delays. Clearly, there are many directions for future research. Indeed, a more comprehensive bifurcation anlysis with regard to the different parameters of the system is a next step and subject of our ongoing work. Given the experimental relevance of the time delays involved in human movement coordination dynamics, it would be of particular importance to study further the effects of heterogeneity in the intrinsic properties and coupling strengths of the two oscillators in the presence of delays. Another interesting future direction to explore would be the analysis of a two-tier model (different type of coupling) \cite{Banerjee2007}, which has some neurological support; see e.g. \cite{Maki2008}.

\section*{Acknowledgement}
This work was funded by the European Project AlterEgo FP7 ICT 2.9 - Cognitive Sciences and Robotics, Grant Number 600610, and KT-A was supported by grant EP/L000296/1 of the Engineering and Physical Sciences Research Council (EPSRC).


\begin{thebibliography}{32}
\providecommand{\natexlab}[1]{#1}
\providecommand{\url}[1]{{#1}}
\providecommand{\urlprefix}{URL }
\expandafter\ifx\csname urlstyle\endcsname\relax
  \providecommand{\doi}[1]{DOI~\discretionary{}{}{}#1}\else
  \providecommand{\doi}{DOI~\discretionary{}{}{}\begingroup
  \urlstyle{rm}\Url}\fi
\providecommand{\eprint}[2][]{\url{#2}}

\bibitem[{Avitabile et~al(2016)Avitabile, Slowinski, Bardy, and
  Tsaneva-Atanasova}]{Avitabile2016}
Avitabile D, Slowinski P, Bardy B, Tsaneva-Atanasova K (2016) Beyond in-phase
  and anti-phase coordination in a model of joint action. Biological
  Cybernetics (to appear)

\bibitem[{Banerjee and Jirsa(2007)}]{Banerjee2007}
Banerjee A, Jirsa VK (2007) How do neural connectivity and time delays
  influence bimanual coordination? Biological Cybernetics 96(2):265--278

\bibitem[{Bourbousson et~al(2010)Bourbousson, Seve, and
  McGarry}]{Bourbousson2010}
Bourbousson J, Seve C, McGarry T (2010) Space--time coordination dynamics in
  basketball: Part 1. intra-and inter-couplings among player dyads. Journal of
  Sports Sciences 28(3):339--347

\bibitem[{Calvin and Jirsa(2011)}]{Calvin2011}
Calvin S, Jirsa VK (2011) Perspectives on the Dynamic Nature of Coupling in
  Human Coordination, Springer-Verlag, pp 91--114

\bibitem[{Diekmann et~al(1995)Diekmann, Van~Gils, Lunel, and
  Walther}]{Diekmann1995}
Diekmann O, Van~Gils SA, Lunel SM, Walther HO (1995) Delay Equations:
  Functional-, Complex-, and Nonlinear Analysis. Springer-Verlag

\bibitem[{Duarte et~al(2012)Duarte, Ara{\'u}jo, Davids, Travassos, Gazimba, and
  Sampaio}]{Duarte2012}
Duarte R, Ara{\'u}jo D, Davids K, Travassos B, Gazimba V, Sampaio J (2012)
  Interpersonal coordination tendencies shape 1-vs-1 sub-phase performance
  outcomes in youth soccer. Journal of Sports Sciences 30(9):871--877

\bibitem[{Fabre-Thorpe et~al(2001)Fabre-Thorpe, Delorme, Marlot, and
  Thorpe}]{FabreThorpe2001}
Fabre-Thorpe M, Delorme A, Marlot C, Thorpe S (2001) A limit to the speed of
  processing in ultra-rapid visual categorization of novel natural scenes. Journal of
  Cognitive Neuroscience 13(2):171--180, \doi{10.1162/089892901564234}

\bibitem[{Haken et~al(1985)Haken, Kelso, and Bunz}]{Haken1985}
Haken H, Kelso JS, Bunz H (1985) A theoretical model of phase transitions in
  human hand movements. Biological Cybernetics 51(5):347--356

\bibitem[{Kay et~al(1987)Kay, Kelso, Saltzman, and Sch\~ner}]{Kay1987}
Kay BA, Kelso JA, Saltzman EL, Sch\~ner G (1987) Space-time behavior of single
  and bimanual rhythmical movements: data and limit cycle model. Journal of
  Experimental Psychology, Human Perception and Performance 13(2):178--92

\bibitem[{Kelso(1997)}]{Kelso1997}
Kelso JS (1997) Dynamic Patterns: the Self-Organization of Brain and Behavior.
  MIT press

\bibitem[{Kelso et~al(2009)Kelso, de~Guzman, Reveley, and Tognoli}]{Kelso2009}
Kelso JS, de~Guzman GC, Reveley C, Tognoli E (2009) Virtual partner interaction
  (vpi): exploring novel behaviors via coordination dynamics. PLoS One
  4(6):e5749

\bibitem[{Kralemann et~al(2008)Kralemann, Cimponeriu, Rosenblum, Pikovsky, and
  Mrowka}]{Kralemann2008}
Kralemann B, Cimponeriu L, Rosenblum M, Pikovsky A, Mrowka R (2008) Phase
  dynamics of coupled oscillators reconstructed from data. Physical Review E
  77(6):066,205

\bibitem[{Li et~al(2012)Li, Levine, and Loeb}]{Li2012}
Li Y, Levine WS, Loeb GE (2012) A two-joint human posture control model with
  realistic neural delays. IEEE
  Transactions on Neural Systems and Rehabilitation Engineering 20(5):738--748

\bibitem[{Loram et~al(2006)Loram, Gawthrop and Lakie}]{Loram2006}
Loram ID, Gawthrop PJ, Lakie M (2006) The frequency of human, manual adjustments
  in balancing an inverted pendulum is constrained by intrinsic physiological
  factors. The Journal of physiology 577(1):417--432

\bibitem[{Maki et~al(2008)Maki, Wong, Sugiura, Ozaki, and Sadato}]{Maki2008}
Maki Y, Wong KFK, Sugiura M, Ozaki T, Sadato N (2008) Asymmetric control
  mechanisms of bimanual coordination: an application of directed connectivity
  analysis to kinematic and functional {MRI} data. Neuroimage 42(4):1295--1304

\bibitem[{Marten et~al(2009)Marten, Rodrigues, Benjamin, Richardson, and
  Terry}]{Marten2009}
Marten F, Rodrigues S, Benjamin O, Richardson MP, Terry JR (2009) Onset of
  polyspike complexes in a mean-field model of human electroencephalography and
  its application to absence epilepsy. Philosophical Transactions of the Royal
  Society of London A: Mathematical, Physical and Engineering Sciences
  367(1891):1145--1161

\bibitem[{Miall et~al(1993)Miall, Weir, and Stein}]{Miall1993}
Miall RC, Weir D, Stein J (1993) Intermittency in human manual tracking tasks.
  Journal of Motor Behavior 25(1):53--63

\bibitem[{Michel et~al(2004)Michel, Seeck, and Murray}]{Michel2004}
Michel C, Seeck M, Murray M (2004) The speed of visual cognition. Supplements
  to Clinical Neurophysiology 57:617--627

\bibitem[{Molenaar and Newell(2003)}]{Molenaar2003}
Molenaar P, Newell KM (2003) Direct fit of a theoretical model of phase
  transition in oscillatory finger motions. British Journal of Mathematical and
  Statistical Psychology 56(2):199--214

\bibitem[{Mortl et~al(2012)Mortl, Lorenz, Vlaskamp, Gusrialdi, Schuba, and
  Hirche}]{Mortl2012}
Mortl A, Lorenz T, Vlaskamp BN, Gusrialdi A, Schuba A, Hirche S (2012) Modeling
  inter-human movement coordination: synchronization governs joint task
  dynamics. Biological Cybernetics 106(4-5):241--259

\bibitem[{Oullier et~al(2008)Oullier, de~Guzman, Jantzen, Lagarde, and
  Kelso}]{Oullier2008}
Oullier O, de~Guzman GC, Jantzen KJ, Lagarde J, Kelso JAS (2008) Social
  coordination dynamics: measuring human bonding. Social Neuroscience
  3(2):178--192

\bibitem[{Peper and Beek(1998)}]{Peper1998}
Peper C, Beek PJ (1998) Distinguishing between the effects of frequency and
  amplitude on interlimb coupling in tapping a 2:3 polyrhythm. Experimental
  Brain Research 118(1):78--92

\bibitem[{Post et~al(2000)Post, Peper, Daffertshofer, and Beek}]{Post2000}
Post A, Peper C, Daffertshofer A, Beek PJ (2000) Relative phase dynamics in
  perturbed interlimb coordination: stability and stochasticity. Biological
  Cybernetics 83(5):443--459

\bibitem[{Schmidt and O'Brien(1997)}]{Schmidt1997}
Schmidt R, O'Brien B (1997) Evaluating the dynamics of unintended interpersonal
  coordination. Ecological Psychology 9(3):189--206

\bibitem[{Sch{\"o}ner et~al(1986)Sch{\"o}ner, Haken, and Kelso}]{Schoner1986}
Sch{\"o}ner G, Haken H, Kelso J (1986) A stochastic theory of phase transitions
  in human hand movement. Biological Cybernetics 53(4):247--257

\bibitem[{Scott~Kelso et~al(1983)Scott~Kelso, Putnam, and Goodman}]{Kelso1983}
Scott~Kelso J, Putnam CA, Goodman D (1983) On the space-time structure of human
  interlimb co-ordination. The Quarterly Journal of Experimental Psychology
  35(2):347--375

\bibitem[{{Sieber} et~al(2014){Sieber}, {Engelborghs}, {Luzyanina}, {Samaey},
  and {Roose}}]{DDEBIFTOOL}
{Sieber} J, {Engelborghs} K, {Luzyanina} T, {Samaey} G, {Roose} D (2014)
  {DDE-BIFTOOL Manual - Bifurcation analysis of delay differential equations}.
  ArXiv e-prints \eprint{1406.7144}

\bibitem[{Tass et~al(1995)Tass, Wunderlin, and Schanz}]{Tass1995}
Tass P, Wunderlin A, Schanz M (1995) A theoretical model of sinusoidal forearm
  tracking with delayed visual feedback. Journal of Biological Physics
  21(2):83--112

\bibitem[{Tass et~al(1996)Tass, Kurths, Rosenblum, Guasti, and
  Hefter}]{Tass1996}
Tass P, Kurths J, Rosenblum M, Guasti G, Hefter H (1996) Delay-induced
  transitions in visually guided movements. Physical Review E 54(3):R2224

\bibitem[{Thorpe et~al(1996)Thorpe, Fize, Marlot et~al}]{Thorpe1996}
Thorpe S, Fize D, Marlot C, et~al (1996) Speed of processing in the human
  visual system. Nature 381(6582):520--522

\bibitem[{Varlet et~al(2012)Varlet, Marin, Raffard, Schmidt, Capdevielle,
  Boulenger, Del-Monte, and Bardy}]{Varlet2012}
Varlet M, Marin L, Raffard S, Schmidt RC, Capdevielle D, Boulenger JP,
  Del-Monte J, Bardy BG (2012) Impairments of social motor coordination in
  schizophrenia. PLoS One 7(1):e29,772

\bibitem[{Varlet et~al(2014)Varlet, Marin, Capdevielle, Del-Monte, Schmidt,
  Salesse, Boulenger, Bardy, and Raffard}]{Varlet2014}
Varlet M, Marin L, Capdevielle D, Del-Monte J, Schmidt R, Salesse R, Boulenger
  JP, Bardy BG, Raffard S (2014) Difficulty leading interpersonal coordination:
  Towards an embodied signature of social anxiety disorder. Frontiers in
  Behavioral Neuroscience 8(29)

\bibitem[{Warren(2006)}]{Warren2006}
Warren WH (2006) The dynamics of perception and action. Psychological Review
  113(2):358--89

\end{thebibliography}


\end{document}